\PassOptionsToPackage{unicode}{hyperref}
\PassOptionsToPackage{hyphens}{url}
\documentclass[
  british,
  11pt,
]{article}
\usepackage{xcolor}
\usepackage{graphicx}
\usepackage[margin=1in]{geometry}
\usepackage{amsmath,amssymb}
\usepackage{iftex}
\ifPDFTeX
  \usepackage[T1]{fontenc}
  \usepackage[utf8]{inputenc}
  \usepackage{textcomp} 
\else 
  \usepackage{unicode-math} 
  \defaultfontfeatures{Scale=MatchLowercase}
  \defaultfontfeatures[\rmfamily]{Ligatures=TeX,Scale=1}
\fi
\usepackage{lmodern}
\ifPDFTeX\else
\fi
\IfFileExists{upquote.sty}{\usepackage{upquote}}{}
\IfFileExists{microtype.sty}{
  \usepackage[]{microtype}
  \UseMicrotypeSet[protrusion]{basicmath} 
}{}
\makeatletter
\@ifundefined{KOMAClassName}{
  \IfFileExists{parskip.sty}{%
    \usepackage{parskip}
  }{
    \setlength{\parindent}{0pt}
    \setlength{\parskip}{6pt plus 2pt minus 1pt}}
}{
  \KOMAoptions{parskip=half}}
\makeatother
\ifLuaTeX
\usepackage[bidi=basic,shorthands=off]{babel}
\else
\usepackage[bidi=default,shorthands=off]{babel}
\fi
\ifLuaTeX
  \usepackage{selnolig} 
\fi
\providecommand{\tightlist}{%
  \setlength{\itemsep}{0pt}\setlength{\parskip}{0pt}}
\usepackage{bookmark}
\IfFileExists{xurl.sty}{\usepackage{xurl}}{} 
\makeatletter
\@ifundefined{xmpquote}{}{}
\makeatother
\hypersetup{
  pdftitle={Exact majority C-colourings of balanced Hamming graphs and grids},
  pdfauthor={ZeHua LU},
  pdflang={en-GB},
  hidelinks,
  pdfcreator={LaTeX via pandoc}}

\title{Exact majority C-colourings of balanced Hamming graphs and grids}
\author{ZeHua LU\\
Independent researcher\\
\texttt{zhouyiqian19@gmail.com}}
\date{15 September 2026}

\begin{document}
\maketitle


\subsection{Abstract}\label{abstract}

A majority C-colouring partitions a graph into classes in which every
vertex has at least half of its neighbours. Write \(M(G)\) for the
maximum number of classes. We determine \[
M(K_q^{\square(2k+1)})=
\left\lfloor\frac{q^{k+1}}{\lfloor q/2\rfloor+1}\right\rfloor
\qquad(q\ge3,\ k\ge0).
\] The lower bound follows from explicit rectangular partitions and a
uniform three-dimensional bridge. A punctured-rectangle construction
supplies both the bridge and a partition in the intermediate dimension.
The upper bound is a classical consequence of Hamming edge isoperimetry;
an elementary proof is included. We also prove
\(M(C_m\square P_n)=\frac n2\lfloor m/2\rfloor\) for \(m\ge4\) and even
\(n\ge2\), contradicting the cylinder assertion of Conjecture 4 in
arXiv:2608.27669v1 for odd \(m\ge7\) and even \(n\ge6\). Further layer
bounds determine additional cylinder and torus families.

\subsection{1. Definitions and main
results}\label{definitions-and-main-results}

We study majority C-colourings as introduced in {[}2{]} and considered
for Cartesian products in {[}1{]}. All graphs are finite and simple. Put
\(d_S(v)=|N_G(v)\cap S|\). A nonempty set \(S\subseteq V(G)\) is a
\textbf{majority class} if \[
d_S(v)\ge\lceil d_G(v)/2\rceil\qquad(v\in S).
\] A majority C-colouring is a vertex partition into majority classes;
\(M(G)\) is its largest possible number of parts. Equivalently, each
part satisfies \(d_S(v)\ge d_{V(G)\smallsetminus S}(v)\) on its own vertices.
This is the condition for a strong defensive alliance {[}11, Sections
1.2 and 1.3{]}, also called a defensive \(0\)-alliance {[}7, Section 1;
8{]}. No domination condition is imposed. Thus \[
M(G)=\psi_0^d(G),
\] the defensive \(0\)-alliance partition number. On graphs whose
degrees are all even, the defensive \(-1\) and \(0\) conditions
coincide. This includes tori and Hamming graphs with odd \(q\), but not
the degree-three boundary vertices of cylinders.

We use \(P_n\) for the path on \(n\) vertices, \(C_m\) for the
\(m\)-cycle, and \(\square\) for Cartesian product. Set \[
H(q,D)=K_q^{\square D},\quad H(q,0)=K_1,\quad s(q)=\lfloor q/2\rfloor+1.
\] An \textbf{axial clique} in a Hamming graph has all but one
coordinate fixed. A partition into induced subgraphs refers to their
vertex sets; it does not remove edges between the parts.

Bujtás, Dettlaff, Furmańczyk and Laskowska {[}1, Open problem 3{]} ask
for exact values in odd-dimensional balanced Hamming graphs. Their
Conjecture 4 asserts, for \(m,n\ge6\) with at least one odd, \[
M(C_m\square P_n)=M(C_m\square C_n)
=1+\lfloor m/2\rfloor\lfloor n/2\rfloor.
\] Throughout, references to these problems concern arXiv:2608.27669v1.

\textbf{Theorem 1 (balanced Hamming graphs).} For every integer
\(q\ge3\) and \(k\ge0\), \[
\boxed{M(H(q,2k+1))=\left\lfloor\frac{q^{k+1}}{s(q)}\right\rfloor.}
\] This answers Open problem 3 of {[}1{]}. The proof is constructive for
every alphabet size and dimension.

\textbf{Theorem 2 (cylinders).} Let \(m\ge4,n\ge2\), and
\(a=\lfloor m/2\rfloor\).

\begin{enumerate}
\def\labelenumi{\arabic{enumi}.}
\tightlist
\item
  If \(n\) is even, then \(M(C_m\square P_n)=an/2\).
\item
  If \(n\) is odd, then \(M(C_m\square P_n)=1+a(n-1)/2\) whenever \(m\)
  is even, \(m\le18\), or \(n\in\{3,5,7,9\}\).
\end{enumerate}

The even-\(m\), odd-\(n\) case is already included in {[}1, Theorem 12
and Remark 13{]}; we give a common short proof. The first assertion
disproves the cylinder equality quoted above. The smallest permitted
example is \[
M(C_7\square P_6)=9\ne10.
\] \textbf{Theorem 3 (tori).} If \(m\ge5\) is odd and \(n\ge4\) is even,
then \[
M(C_m\square C_n)=1+\frac n2\lfloor m/2\rfloor.
\] The transposed statement holds. In the odd-by-odd range of Conjecture
4, its formula also holds for \(C_7\square C_n\) for every odd
\(n\ge7\), and for \((m,n)=(9,9),(9,11),(9,13)\) and their transposes.

The main constructive tool is an explicit partition of a rectangle into
axial cliques. In the Hamming proof, one exceptional clique carries the
remainder in the vertex count. A punctured rectangle and a
three-dimensional bridge handle the only remainder that prevents a
one-coordinate step. Sections 2 to 4 prove Theorem 1. Sections 5 and 6
develop additive layer bounds for cylinders and tori; the two appendices
treat the remaining stated small-width cases.

The underlying alliance parameter and several ingredients are
established. The value \(M(K_q\square K_q)=q\) appears in {[}7, Example
1 and the equality discussion following the basic bounds{]}, and the
product rule for defensive alliances is {[}7, Theorem 12{]}. The Hamming
size bound in Section 2 follows from classical edge isoperimetry {[}4,
5, 6{]}. As recalled in {[}1{]}, the binary hypercube case was treated
in {[}3{]}. Our Hamming argument supplies attaining partitions for all
\(q\ge3\).

The grid problem also has a longer history. Lachniet {[}11, Section 4.1,
Problems 5 and 6{]} asked for the ordinary alliance partition numbers of
cylinders and tori. Theorem 2 determines specified families for the
strong-alliance analogue of the cylinder problem. Since tori have even
degree, Theorem 3 gives a partial answer to the same parameter in his
torus problem. Ordinary alliance partitions of path grids were studied
in {[}10{]} and {[}11, Section 2.3{]}. The use of cycles to bound
alliance sizes occurs in {[}9, Proposition 3.4{]}, and cycle and
row-column counting appears in {[}11, Lemma 2.13 and the proof of
Theorem 2.14{]}. Our layer inequalities give an additive refinement
adapted to the present products.

\subsection{2. A classical Hamming size
bound}\label{a-classical-hamming-size-bound}

For integers \(r\le0\) put \(f_q(r)=1\). For \(r\ge0\) write \[
r=a(q-1)+b,\quad 0\le b<q-1,\qquad f_q(r)=q^a(b+1).
\] The definitions agree at zero.

\textbf{Lemma 4.} Let \(q\ge3\) and \(D\ge0\) be integers, and let
\(r\in\mathbb Z\). A nonempty induced subgraph of \(H(q,D)\) with
minimum degree at least \(r\) has at least \(f_q(r)\) vertices.

The lemma follows from Lindsey's classical edge-isoperimetric theorem
{[}4{]}, as presented in {[}5, Theorem 2.3{]}, together with the
standard entropy bound discussed in {[}6, Section 3.1{]}. We include an
elementary proof of the minimum-degree assertion.

\emph{Proof.} Induct on \(D\), simultaneously for all integer
thresholds. Nonpositive thresholds are immediate. In dimension zero a
positive threshold is impossible. Slice the set \(S\) along one
coordinate and let \(t\) be the number of nonempty slices,
\(1\le t\le q\). Each vertex has at most one neighbour in each other
slice. Hence each slice has minimum degree at least \(r-t+1\), and \[
|S|\ge t f_q(r-t+1).
\] For \(r=a(q-1)+b\ge0\) put \(h=t-1\). If \(h\le b\), division by
\(q^a\) reduces the required inequality to \[
(h+1)(b-h+1)-(b+1)=h(b-h)\ge0.
\] If \(h>b\) and \(a\ge1\), division by \(q^{a-1}\) gives \[
t(q+b-t+1)-q(b+1)=(t-b-1)(q-t)\ge0.
\] If \(h>b\) and \(a=0\), the new threshold is negative and
\(t\ge b+2>b+1\). This proves the induction. \(\square\)

The majority threshold in \(H(q,2k+1)\) is \[
\left\lceil(2k+1)(q-1)/2\right\rceil=k(q-1)+s(q)-1.
\] For \(q\ge3\), every class therefore has at least \(q^k s(q)\)
vertices. Division into the \(q^{2k+1}\) vertices gives the upper bound
of Theorem 1.

\subsection{3. Rectangular partitions and the residue
criterion}\label{rectangular-partitions-and-the-residue-criterion}

\textbf{Lemma 5.} Let \(s,a,q\) be positive integers with
\(s\le a\le q<2s\) and write \(qa=Ns+e\), \(0\le e<s\). If \(e\le q-s\),
then \(K_q\square K_a\) has a vertex partition into \(N\) axial cliques.
All have size \(s\), except for exactly one of size \(s+e\) when
\(e>0\).

\emph{Proof.} Use \(q\) rows and \(a\) columns and put \[
h=N-a,\qquad d=s-q+h.
\] We have \(q-s\le h<q\) and \(d\ge0\). If \(h=0\), then \(q=s=a\) and
the columns give the required partition. Assume \(h>0\). Designate \(h\)
rows as special and the other \(q-h\) as common. On the special rows
choose a zero-one matrix with column sums \[
c_0=d+e,\qquad c_1=\cdots=c_{a-1}=d.
\] The assumptions give \(0\le c_j\le h\), and \[
\sum_jc_j=ad+e=h(a-s).
\] Assign each column's marks consecutively, with row indices continuing
cyclically modulo \(h\) from one column to the next. A column never
repeats a row, and each row receives exactly \(a-s\) marks.

In each special row the unmarked entries form a horizontal \(K_s\). In
each column combine all common entries and marked special entries into a
vertical clique. The latter sizes are \(q-h+d=s\), except for column
zero, of size \(s+e\). Every entry is used once; there are \(h+a=N\)
parts. \(\square\)

\textbf{Proposition 6.} Let \(q\ge3\) and \(L\ge1\) be integers, and put
\(s=s(q)\). If \(q^j\bmod s\le q-s\) for \(1\le j\le L\), then
\(H(q,L)\) has a partition into \(\lfloor q^L/s\rfloor\) axial cliques
of size \(s\), except possibly one of size \(s+(q^L\bmod s)\).

\emph{Proof.} For \(L=1\) use the whole \(K_q\). When passing from
dimension \(j\) to \(j+1\), copy every ordinary \(K_s\) for the \(q\)
new coordinate values. Replace the exceptional \(K_a\square K_q\), if
present, using Lemma 5. Its remainder is \(qa\bmod s=q^{j+1}\bmod s\).
If the remainder is zero there is no further exception. The parts remain
axial, and their sizes give the stated count. \(\square\)

For \(L=k+1\) lift each part \(B\) to \(V(H(q,k))\times B\). Its
internal minimum degree is at least \(k(q-1)+s-1\), proving Theorem 1
under the residue condition. This is an instance of the alliance product
rule {[}7, Theorem 12{]}.

For odd \(q=2s-1\), all remainders are allowed. For \(q\equiv2\pmod4\),
\(s\) is even and all remainders are even, so they avoid the only
forbidden remainder \(s-1\). The case \(q=4\) has constant remainder
\(1\) modulo \(3\). The remaining alphabets are covered by the
construction below.

\subsection{4. Punctured rectangles and a uniform
bridge}\label{punctured-rectangles-and-a-uniform-bridge}

\subsubsection{4.1. One construction for punctures and terminal
blocks}\label{one-construction-for-punctures-and-terminal-blocks}

\textbf{Lemma 7 (punctured rectangle and terminal partition).} Let
\(t\ge2\) be an integer and put \(s=2t+1\), \(q=4t=2s-2\), and
\(a=3t+2=(3s+1)/2\). In a \(q\)-row, \(a\)-column array, delete any
\(2t=s-1\) entries from column zero. The remaining graph partitions into
\(6t\) axial copies of \(K_s\). Consequently the full \(K_a\square K_q\)
partitions into \(6t-1=3s-4\) copies of \(K_s\) and one block consisting
of two copies of \(K_s\) sharing exactly one vertex.

\emph{Proof.} Order all \(2t\) rows without a deleted entry first,
followed by any \(t-1\) rows with a deleted entry. These form \(h=3t-1\)
special rows. The other \(t+1\) rows are common, and all have a deleted
entry in column zero.

Place \(t\) marks in each positive column \(1,\ldots,a-1\), assigning
row indices consecutively and cyclically through the ordered special
rows. Since \(t\le h\), a column has no repeated row. The total number
of marks is \[
t(a-1)=t(3t+1)=th+2t.
\] Thus each of the first \(2t\) special rows has \(t+1\) marks, and
each remaining special row has \(t\).

Take the unmarked, undeleted entries of each special row as a horizontal
clique. A row without a hole contributes \(a-(t+1)=s\) entries; a row
with a hole contributes \(a-1-t=s\). In each positive column combine all
common entries and marked special entries as a vertical clique of size
\((t+1)+t=s\). Every undeleted entry is covered once, including column
zero, whose remaining entries all lie in special rows without a hole.
There are \(h+(a-1)=6t\) cliques.

For the terminal assertion, choose a horizontal clique in any row
without a hole, and add all the deleted entries to its block. The
horizontal clique meets column zero at one vertex. That vertex together
with the restored \(s-1\) entries forms a vertical \(K_s\). The union
consists of two \(K_s\)'s sharing one vertex, has size \(2s-1\), and has
minimum degree \(s-1\). All other blocks remain unchanged. \(\square\)

\subsubsection{4.2. A bridge for every remaining
alphabet}\label{a-bridge-for-every-remaining-alphabet}

\textbf{Lemma 8 (uniform bridge).} With the parameters of Lemma 7, \[
K_a\square K_q\square K_q
\] has a partition into \(24t^2+4t-3\) axial copies of \(K_s\) and one
axial \(K_{s+2}\).

\emph{Proof.} View the product as \(q\) layers of the array in Lemma 7.
Fix \(q-2\) positions in column zero as possible holes. Require one
position to be a hole in \(s+2\) layers, and each of the other \(q-3\)
positions to be a hole in \(s\) layers.

Construct these incidence sets by assigning layer indices successively
and cyclically modulo \(q\), continuing from one position to the next.
No position receives the same index twice, because \(s+2\le q\) for
\(t\ge2\). Furthermore, \[
(s+2)+(q-3)s=(q-2)s+2=q(s-1).
\] Every layer therefore contains exactly \(s-1\) holes. Apply Lemma 7
separately to each layer, obtaining \(6tq\) ordinary cliques. Along the
third coordinate the holes themselves form \(q-3\) further \(K_s\)'s and
one \(K_{s+2}\). These cover precisely the entries omitted from the
layer partitions. The number of ordinary cliques is \[
6tq+q-3=24t^2+4t-3,
\] and all are axial. \(\square\)

An independent short construction applies when \(t\ge3\). Lemma 5 first
partitions \(K_q\square K_q\) into \((4s-9)K_s\) and \(K_{s+4}\),
because \(q^2\equiv4\pmod s\) and \(4\le q-s\). Copy the ordinary blocks
along \(K_a\). The remaining rectangle \(K_{s+4}\square K_a\) again
satisfies Lemma 5: \[
s\le s+4\le a<2s,\qquad
a(s+4)\equiv2\pmod s,\qquad 2\le a-s.
\] It leaves one \(K_{s+2}\). The hole construction above also includes
the boundary case \(t=2\).

\subsubsection{4.3. Completing all alphabets and
dimensions}\label{completing-all-alphabets-and-dimensions}

It remains to supplement Proposition 6 when \(q\equiv0\pmod4\). The case
\(q=4\) has constant remainder \(1\) modulo \(3\), and is already
settled. Assume \(q\ge8\), so that \(s=q/2+1\) is odd and
\(q\equiv-2\pmod s\).

Maintain a partition of \(H(q,L)\) into ordinary axial \(K_s\)'s and one
exceptional axial clique of size \(s+e_L\), where \(e_L=q^L\bmod s\).
These residues are nonzero because \(\gcd(q,s)=1\). The initial state
\(L=1\) is the whole \(K_q\). If the next residue is at most
\(q-s=s-2\), add one coordinate using Lemma 5, copying every ordinary
block for each value of the new coordinate.

If the next residue is \(s-1\), the current residue is uniquely
\((s+1)/2\), since \(-2e_L\equiv-1\pmod s\). Hence the current
exceptional clique has size \(a=(3s+1)/2\). From this same state
construct two outputs:

\begin{enumerate}
\def\labelenumi{\arabic{enumi}.}
\tightlist
\item
  Add one coordinate and apply the terminal part of Lemma 7. This gives
  the required partition in dimension \(L+1\), with a single block of
  size \(2s-1\) and minimum degree \(s-1\).
\item
  For continued recursion, instead add two coordinates and apply Lemma
  8. The exceptional block in dimension \(L+2\) is again an axial
  clique, now of size \(s+2\).
\end{enumerate}

The residue at \(L+2\) is \(q(s-1)\equiv2\pmod s\), and \(s+2\le q\), so
the second output restores the invariant. Ordinary blocks are copied
\(q\) or \(q^2\) times as appropriate. Repeating these operations
increases the dimension strictly, and the terminal outputs cover every
skipped dimension. No nonclique terminal block is ever used as input to
the axial recurrence.

Every output has \(q^L\) vertices, ordinary blocks of size \(s\), and an
exceptional block of size \(s+e_L\) if the remainder is nonzero. The
number of blocks is thus \(\lfloor q^L/s\rfloor\), and every block has
minimum degree at least \(s-1\). A zero residue in the other alphabet
cases simply removes the exception, after which all ordinary blocks can
be copied.

For \(L=k+1\), lift each block by \(H(q,k)\) as in Section 3. Its
internal minimum degree is at least \[
k(q-1)+s-1=\left\lceil(2k+1)(q-1)/2\right\rceil.
\] Together with the classical upper bound, this proves Theorem 1 for
every \(q\ge3\).

An exceptional axial clique need not use consecutive coordinate values.
Label its values by any bijection with \(\{0,\ldots,a-1\}\) before
applying a rectangle or bridge, and map back afterwards. This preserves
axiality in the full Hamming graph.

\subsection{5. Layer bounds for
cylinders}\label{layer-bounds-for-cylinders}

\subsubsection{5.1. Even heights and the
counterexample}\label{even-heights-and-the-counterexample}

Write a vertex of \(C_m\square P_n\) as \((x,y)\) with
\(x\in\mathbb Z_m\) and \(0\le y<n\). The \(C_m\) layer with fixed \(y\)
is called a row, denoted \(L_y\). Ambient degrees are three or four, so
every majority class has induced minimum degree at least two.

\textbf{Lemma 9 (row budget).} For \(m\ge4,n\ge2\) and any majority
class \(S\), \[
B(S):=\sum_y\left\lfloor |S\cap L_y|/2\right\rfloor\ge2.
\]

\emph{Proof.} The class contains a simple cycle. If it occupies a single
row, it is the whole \(C_m\) and contributes at least two. Otherwise a
vertex in either extreme occupied row of that cycle has at most one
cycle neighbour outside its row. Both extreme rows therefore contain a
horizontal cycle edge and contribute at least one. \(\square\)

For a partition into \(c\) classes, \[
2c\le\sum_S B(S)\le n\lfloor m/2\rfloor.
\] When \(n\) is even, pair the rows. Within each pair, partition the
columns into intervals of size two, with one interval of size three if
\(m\) is odd. Each resulting \(2\times2\) or \(2\times3\) block has
minimum degree at least two. Giving all blocks different colours attains
\(\frac n2\lfloor m/2\rfloor\), proving Theorem 2(1).

For \(C_7\square P_6\), an optimal nine-colouring is \[
\begin{matrix}
A&A&B&B&C&C&C\\
A&A&B&B&C&C&C\\
D&D&E&E&F&F&F\\
D&D&E&E&F&F&F\\
G&G&H&H&I&I&I\\
G&G&H&H&I&I&I
\end{matrix}
\] with rows wrapping horizontally and the vertical direction a path.
The budget inequality gives \(2c\le18\), ruling out the ten colours
asserted by {[}1, Conjecture 4{]}.

This does not refute the torus equality. Closing the vertical direction
changes the boundary condition: one complete column and nine squares
give ten colours on \(C_7\square C_6\).

\subsubsection{5.2. Odd heights}\label{odd-heights}

Let \(n=2t+1\) and \(a=\lfloor m/2\rfloor\). One full row together with
paired remaining rows gives \(at+1\) colours. For the reverse bound we
may assume classes are connected: splitting a disconnected class into
its components preserves the majority condition and increases the number
of colours.

If a class contains a full row, its budget is at least \(a\). Hence \[
2(c-1)+a\le(2t+1)a,\qquad c\le at+1.
\] It remains to consider connected classes with no full row. Such a
class occupies an interval of at least two rows. Its bottom row is a
union of paths with at least two vertices each: an isolated vertex would
have at most one same-colour neighbour. Both endpoints of every path
require inward vertical neighbours. Thus the bottom row and its inward
neighbouring row each contain at least two vertices. The same holds at
the top.

Every such class consequently has at least two vertices in each row
parity. The \(t\) odd-indexed rows give \[
c\le\lfloor mt/2\rfloor.
\] This is at most \(at+1\) when \(m\) is even or \(t\le3\).

\textbf{Lemma 10 (alternating imbalance).} For a connected majority
class without a full row, define \[
A(S)=\sum_y(-1)^y|S\cap L_y|.
\] Then \(|A(S)|\le2B(S)-2\). If \(B(S)=2\), then \(A(S)=0\). In
particular \[
|A(S)|\le4(B(S)-2).
\]

\emph{Proof.} Give either row parity the negative sign. Let \(E\) be the
number of vertices with that sign, and let \(o\) be the number of rows
with odd cardinality in \(S\). Since \(A=2B+o-2E\), it suffices to show
\(o\le2E-2\). Write \(\ell\) for the number of occupied rows.

For \(\ell=2\), we have \(E\ge2\) and \(o\le2\). For \(\ell=3\), if the
end rows have negative sign then \(E\ge4\). Otherwise the middle row has
negative sign and \(E\ge2\); if \(E=2\) it is even, so \(o\le2\), and if
\(E\ge3\) the bound is immediate. For \(\ell\ge4\), the first two and
last two rows are disjoint pairs, each containing a negative row with at
least two vertices. Every other occupied negative row contains at least
one vertex. Therefore \(E\ge\lfloor\ell/2\rfloor+2\), giving
\(2E-2\ge\ell\ge o\). Interchanging the signs proves the absolute-value
assertion.

If \(B=2\), the boundary observation rules out three or more occupied
rows. The two row sizes are each two or three. If unequal, the two
vertices in the smaller row must be adjacent and must also be the
endpoints of a three-vertex path in the other row. Those endpoints are
not adjacent in \(C_m\) for \(m\ge4\), a contradiction. Thus the row
sizes agree and \(A=0\). For \(B\ge3\) use \(2B-2\le4(B-2)\).
\(\square\)

Summing \(\sum_S A(S)=m\) gives \[
\left\lceil m/4\right\rceil
\le\sum_S(B(S)-2)
\le(2t+1)a-2c.
\] Combining all cases, we obtain \[
M(C_m\square P_{2t+1})\le
\max\left\{at+1,\,
\min\left(
\left\lfloor mt/2\right\rfloor,
\left\lfloor\frac{(2t+1)a-\lceil m/4\rceil}{2}\right\rfloor
\right)\right\}.
\] For odd \(m=2a+1\le17\), \[
\left\lfloor\frac{a-\lceil(2a+1)/4\rceil}{2}\right\rfloor\le1.
\] Hence \(c\le at+1\). This completes Theorem 2 for \(m\le18\).
Appendix A supplies the remaining height-nine case.

The factor four in the last inequality of Lemma 10 cannot be replaced by
two without additional assumptions. The perimeter of a \(3\times3\)
rectangle has row sizes \(3,2,3\), giving \(B=3\) and \(|A|=4\).

\subsection{6. Layer bounds for tori}\label{layer-bounds-for-tori}

\subsubsection{6.1. Mixed parities}\label{mixed-parities}

Throughout this section assume \(m,n\ge4\). In \(C_m\square C_n\) call
the \(C_m\) layers rows. If a majority class misses a row, cut the torus
there. Its induced graph does not change, and the cycle argument of
Lemma 9 gives budget at least two.

Suppose a colouring has \(c\) classes, of which \(u\) meet every row.
The remaining classes occupy at most \(m-u\) vertices in each row. Thus
\[
2(c-u)\le n\left\lfloor\frac{m-u}{2}\right\rfloor.
\] Different classes may use different cuts; their budgets are still
counted in the same original rows.

Put \(m=2a+1\), \(n=2b\), with \(a,b\ge2\). According as \(u=2v\) or
\(u=2v+1\), we get \[
c\le ab+v(2-b)\le ab,\qquad
c\le ab+1+v(2-b)\le ab+1.
\] A complete column and a square partition of the remaining \(2a\) by
\(2b\) array attain \(ab+1\), proving the mixed-parity statement.

For two odd lengths, a single colour on the union of a complete row and
a complete column, and separate colours on squares in the remaining
even-by-even array, gives \[
M(C_m\square C_n)\ge1+\lfloor m/2\rfloor\lfloor n/2\rfloor.
\] The preceding upper-bound argument gives \[
M(C_m\square C_n)\le
\left\lfloor\frac{n\lfloor m/2\rfloor}{2}\right\rfloor+(m\bmod2)
\] and the transposed bound. To justify the displayed maximum, set
\(F(u)=u+\lfloor n\lfloor(m-u)/2\rfloor/2\rfloor\). For \(n\ge4\),
\(F(u+2)\le F(u)\), so only \(u=0,1\) need comparison. These bounds do
not settle every odd-by-odd pair.

\subsubsection{6.2. Small-budget classes}\label{small-budget-classes}

Define the row and column budgets of a class by \[
R(S)=\sum_y\lfloor |S\cap L_y|/2\rfloor,\qquad
C(S)=\sum_x\lfloor |S\cap\mathrm{column}_x|/2\rfloor.
\]

\textbf{Lemma 11.} For \(m,n\ge6\), a majority class with \(R<2\) is
exactly a complete straight column and has \(R=0\). The transposed
assertion holds.

\emph{Proof.} Every row contains at most one vertex, except possibly one
row containing two or three. Choose a cycle in the class and remove that
exceptional row, or any row of the cycle if no exceptional row exists.
The cycle cannot lie wholly in the removed row, since fewer than \(m\)
vertices are available there.

Outside the removed row, every nontrivial portion of the cycle is a
vertical path in one column. To return to the removed row, such a path
must traverse all \(n-1\) remaining rows. There can be only one such
path: two would put two vertices in each remaining row. Its two
end-neighbours in the removed row have the same column. The cycle
therefore closes through that one vertex; a nontrivial horizontal return
to the same column would require the entire horizontal cycle. Thus the
selected cycle is a straight complete column.

Extra vertices could occur only in the exceptional row. They have no
vertical neighbours in the class. The at most three occupied vertices of
that row cannot support all extra vertices with degree at least two for
\(m\ge6\). Thus there are none. \(\square\)

For \(m=7\) and odd \(n\ge7\), deleting a straight column class leaves
\(C_n\square P_6\), with at most \(3\lfloor n/2\rfloor\) other classes.
If no such class exists, Lemma 11 gives \(R\ge2\) for every class and
hence \[
M\le\lfloor3n/2\rfloor=1+3\lfloor n/2\rfloor.
\] The cross-and-squares construction attains this value.

\textbf{Lemma 12.} For \(m,n\ge9\), a majority class with \(R=C=2\) is
exactly a unit square.

\emph{Proof.} Each row slice is a proper subset of its cycle and has at
most five vertices. A nonempty proper slice with \(z\) vertices induces
at most \(z-1\le2\lfloor z/2\rfloor\) edges. Therefore a cycle in the
class uses at most four horizontal and four vertical edges. It has
length at most eight, cannot wind around either factor, and lifts to a
simple cycle in the planar square grid. Its coordinate spans are smaller
than the respective periods, so distinct lifted rows or columns do not
merge when computing the original budgets.

Every occupied row and column of that planar cycle contains at least two
vertices. In an extreme row a horizontal edge is necessary. A unique
vertex in an intermediate row would disconnect the rest of the cycle
into nonempty portions above and below it, which is impossible. The same
applies to columns. The budgets thus force two rows and two columns,
giving a unit square.

Every cycle in the class is such a square. Any different unit square,
whether sharing an edge, one vertex, or no vertices, increases at least
one of the two budgets when added to the first. Hence there is only one
cycle. A connected finite unicyclic graph of minimum degree at least two
is its cycle; any additional component would contain another cycle. The
whole class is the square. \(\square\)

On \(9\times9\), a straight row or column class can be deleted, leaving
at most sixteen other classes. Otherwise eighteen classes would force
\(R=C=2\) for all classes, because both total budgets are at most
thirty-six. Lemma 12 would cover only \(18\cdot4=72\) vertices instead
of eighty-one. Thus \(M=17\). Appendix B proves \(M=21\) and \(25\) on
\(9\times11\) and \(9\times13\).

\subsection{7. Further questions and
verification}\label{further-questions-and-verification}

Theorem 1 determines all balanced Hamming graphs in the odd-dimensional
range of {[}1, Open problem 3{]}. It would be useful to find similarly
explicit optimal partitions for products of complete graphs of unequal
orders, or for other required proportions of same-coloured neighbours.

The remaining cylinder and torus cases require additional structural
information. The bounds proved here leave \[
46\le M(C_{19}\square P_{11})\le47,
\] and do not determine every odd-by-odd torus. The cycle and layer
methods of Sections 5 and 6 provide constraints for these cases.

The supplementary material implements the rectangular,
punctured-rectangle and bridge constructions using only Python's
standard library. Its verifier checks complete coverage, block sizes,
axiality, majority degrees after lifting, and exact recurrence counts.
The finite ranges and output records are documented with the code. All
theorems in this paper follow from the mathematical proofs above; the
programs provide reproducible implementation checks.

\textbf{Author contributions and use of AI in the research process.} The
author selected the research problem, set the research objectives, and
decided which approaches to pursue or discard. The author worked through
some of the proofs by hand and checked individual proof arguments step
by step. Larger computational checks were carried out primarily with AI
assistance. OpenAI Codex provided substantial assistance with exploring
constructions, developing proof arguments, searching the literature, and
generating and debugging verification code. The author personally
reviewed and edited the full manuscript as needed. The accompanying
proofs are self-contained; the finite computational checks and
model-based audits do not constitute formal proof verification.

\subsection{Appendix A. Height-nine
cylinders}\label{appendix-a.-height-nine-cylinders}

We prove \(M(C_m\square P_9)=1+4\lfloor m/2\rfloor\). Even \(m\) and
classes containing a full row are already handled. Suppose
\(m=2a+1\ge5\) and all classes are connected without full rows. Give
rows \(0,2,4,6,8\) positive sign and the other four rows negative sign.

If \(c\ge4a+2\), the total of \(4m=8a+4\) negative vertices and the
minimum of two per class force \[
c=4a+2,\qquad E(S)=2\quad\text{for each class}.
\] The boundary observation from Section 5 implies that such a class
meets a single negative row and one or both adjacent positive rows. Four
occupied rows give two distinct negative rows with at least two vertices
each. Three occupied rows with negative end rows likewise give at least
four negative vertices.

Let \(U\) be the two occupied columns in the negative row. Each positive
boundary row is a proper subset of \(C_m\), and every path component has
its endpoints vertically adjacent to \(U\). There can only be one
component, an arc joining the two columns in \(U\).

If those columns are adjacent, the other arc uses a complete row and is
forbidden. Each positive row has two vertices, making the positive count
\(O(S)\) even. If they are not adjacent, both negative vertices require
same-colour neighbours above and below. The two positive arcs are either
the same arc or complementary arcs of \(C_m\). Since \(m\) is odd, an
odd \(O(S)\) forces complementary arcs with total vertex count \(m+2\),
giving \[
B(S)=\lfloor c_{\mathrm{top}}/2\rfloor+1+\lfloor c_{\mathrm{bottom}}/2\rfloor=a+2.
\] The total positive count is \(5m\), which is odd. Some class has this
larger budget. Every other class has budget at least two, so \[
2(c-1)+(a+2)\le9a,\qquad c\le4a,
\] contrary to \(c=4a+2\). The full-row and paired-row construction
attains \(4a+1\).

\subsection{\texorpdfstring{Appendix B. The \(9\times11\) and
\(9\times13\)
tori}{Appendix B. The 9 x 11 and 9 x 13 tori}}\label{appendix-b.-the-9times11-and-9times13-tori}

Let \(n=2b+1\), with \(b=5\) or \(6\). Deleting a straight row or column
class gives an even-height cylinder and at most \(4b+1\) total classes.
Otherwise each class has \(R,C\ge2\). The row budget gives \(M\le4b+2\).

Suppose \(K=4b+2\) classes exist. Then every class has \(R=2\), and \[
\sum_S(C(S)-2)\le9b-2K=b-4=:E.
\] Classes with \(C=2\) are squares. Let \(z\) be the number of
exceptions and \(e\) the sum of their excesses \(C-2\), so
\(z\le e\le E\). Every exception has \(C\le2+E\le4\), hence at most nine
vertices per column and at most five per row. All slices are proper
subsets of their cycles. The slice edge bound of Lemma 12 gives \[
|S|\le |E(G[S])|\le2R(S)+2C(S)=8+2(C(S)-2),
\] where the first inequality follows from minimum degree two. Counting
the full vertex set yields \[
9(2b+1)\le4(K-z)+8z+2e
\le16b+8+6(b-4).
\] For \(b=5\), this says \(99\le94\); for \(b=6\), it says
\(117\le116\). Both are contradictions. The cross-and-squares
construction attains \(4b+1\) in each case.

\subsection{Declaration of generative AI and AI-assisted technologies in
the manuscript preparation
process}\label{declaration-of-generative-ai-and-ai-assisted-technologies-in-the-manuscript-preparation-process}

During the preparation of this work, the author used OpenAI Codex to
assist with mathematical exploration and the development of
constructions and proof arguments, literature searches, the generation
and debugging of verification code, and the drafting and revision of the
manuscript. The author selected the research problem, directed the
investigation, worked through some of the proofs by hand, and checked
individual proof arguments step by step. Larger computational checks
were carried out primarily with AI assistance. After using this tool,
the author reviewed and edited the content as needed and takes full
responsibility for the content of the published article.

\subsection{References}\label{references}

\begin{enumerate}
\def\labelenumi{\arabic{enumi}.}
\tightlist
\item
  Cs. Bujtás, M. Dettlaff, H. Furmańczyk and A. Laskowska,
  \emph{Majority C-coloring in Cartesian products}, arXiv:2608.27669v1
  (2026). \url{https://arxiv.org/html/2608.27669v1}.
\item
  Cs. Bujtás, M. Dettlaff, H. Furmańczyk and A. Laskowska,
  \emph{Majority C-coloring of graphs}, arXiv:2604.20752v1 (2026).
  \url{https://arxiv.org/html/2604.20752v1}.
\item
  Cs. Bujtás, E. Sampathkumar, Zs. Tuza, L. Pushpalatha and R. C.
  Vasundhara, \emph{Improper C-colorings of graphs}, Discrete Applied
  Mathematics \textbf{159} (2011), 174--186.
  \url{https://doi.org/10.1016/j.dam.2010.11.004}.
\item
  J. H. Lindsey II, \emph{Assignment of numbers to vertices}, American
  Mathematical Monthly \textbf{71} (1964), 508--516.
\item
  S. L. Bezrukov, \emph{Edge Isoperimetric Problems on Graphs} (1999),
  author-hosted survey, Theorem 2.3.
  \url{http://cs2.uwsuper.edu/sb/Papers/eip.pdf}.
\item
  S. Diskin and W. Samotij, \emph{Isoperimetry in product graphs},
  arXiv:2407.02058v1 (2024), Section 3.1 and footnote 2.
  \url{https://arxiv.org/abs/2407.02058v1}.
\item
  I. G. Yero, S. Bermudo, J. A. Rodríguez-Velázquez and J. M. Sigarreta,
  \emph{Partitioning a graph into defensive k-alliances}, Acta
  Mathematica Sinica (English Series) \textbf{27} (2011), 73--82.
  Consulted version: arXiv:0901.4923v2.
  \url{https://arxiv.org/abs/0901.4923v2}.
\item
  I. González Yero and J. A. Rodríguez-Velázquez, \emph{Defensive
  alliances in graphs: a survey}, arXiv:1308.2096v1 (2013).
  \url{https://arxiv.org/abs/1308.2096v1}.
\item
  L. Eroh and R. Gera, \emph{Alliance Partition Number in Graphs}, Ars
  Combinatoria \textbf{103} (2012), 519--529. Consulted public
  manuscript: \url{https://hdl.handle.net/10945/41332}.
\item
  T. W. Haynes and J. A. Lachniet, \emph{The alliance partition number
  of grid graphs}, AKCE International Journal of Graphs and
  Combinatorics \textbf{4} (2007), 51--59.
\item
  J. Lachniet, \emph{Alliance Partitions in Graphs}, M.S. thesis, East
  Tennessee State University (2007). Sections 1.2 and 1.3, 2.3 and 4.1.
  \url{https://dc.etsu.edu/etd/2080/}.
\end{enumerate}

\end{document}